\magnification=\magstephalf
\input amstex
\loadbold
\documentstyle{amsppt}
\refstyle{A}
\NoBlackBoxes
\vsize=7.5in
\def\pf{\hfill $\square$}
\def\c{\cite}

\def\fg{\frak{g}}

\def\fk{\frak{k}}

\def\end{\text{End}}

\def\TF_1{\widetilde F_{1}}
\def\TF_2{\widetilde F_{2}}
\def\TF{\widetilde F}
\def\TG{\widetilde G}
\def\tg{\widetilde \fg}

\topmatter
\title The complete integrability of a Lie-Poisson system proposed
by Bloch and Iserles  \endtitle
\leftheadtext{L.-C. Li and C. Tomei}
\rightheadtext{Lie-Poisson system}
\author Luen-Chau Li and Carlos Tomei\endauthor
\address{L.-C. Li, Department of Mathematics,Pennsylvania State University
University Park, PA  16802, USA}\endaddress
\email luenli\@math.psu.edu\endemail
\address{C. Tomei, Departamento de Matem\'atica, PUC-Rio, R. Marqu\^es
de S. Vicente 225, Rio de Janeiro, RJ 22453-900, Brazil}\endaddress
\email tomei\@mat.puc-rio.br\endemail \dedicatory{Dedicated with
affect and admiration to Percy Deift on the occasion of his 60th
birthday}\enddedicatory \abstract  We establish the Liouville
integrability of the differential equation  $\dot S(t)= [N,S^2(t)],$
recently considered by Bloch and Iserles. Here, $N$ is a real,
fixed, skew-symmetric matrix and $S$ is real symmetric. The equation
is realized as a Hamiltonian vector field on a coadjoint orbit of a
loop group, and sufficiently many commuting integrals are presented,
together with a solution formula for their related flows in terms of
a Riemann-Hilbert factorization problem. We also answer a question
raised by Bloch and Iserles, by realizing the same system on a
coadjoint orbit of a finite dimensional Lie group.

\endabstract
\endtopmatter

\nobreak{\bf Keywords: } Integrable systems, Lax pairs, Toda flows,
loop groups. \hfill\break

\nobreak{\bf MSC-class: } 37K10, 70H06; 35Q58, 22E67 \hfill

\document
\subhead
1. \ Introduction
\endsubhead
\bigskip

In a recent paper \c{BI}, Bloch and Iserles considered the differential
equation
$$\dot S(t)= [N,S^2(t)],  \eqno{(*)}$$
where both $S$ and $N$ are real $n \times n$ matrices, $S$ is
symmetric and $N$ is a constant skew-symmetric matrix. The
similarity with the celebrated Toda differential equation (\c{F},
\c{Mo}) is hard to miss, and one might expect that both equations
have some common properties. Indeed, from the equivalent Lax pair
formulation
$$\dot S = [NS+SN, S] \quad $$
Bloch and Iserles deduced that $S$ evolves by orthogonal
conjugation: this implies global existence of solutions and spectrum
invariance of $S(t)$. Also, following Manakov's approach to the
rigid body equation \c{Ma}, they showed that (*) is equivalent to
the following Lax equation with spectral parameter
$$(S + z N )^{.} = [ NS + SN + z N^2 , S + z N],$$
from which a large collection of conserved quantities is clear, namely,
the nontrivial coefficients in $z$ in the expansions of
$tr (S + z N)^{k}$. Even more, they showed that equation (*) is a
Lie-Poisson system on an appropriate phase space and, corroborated
by numerical experiments, conjectured the complete integrability of
the system.

In this paper, we shall formulate equation (*) as a Hamiltonian system on a
coadjoint orbit of a loop group, in the spirit of \c{DLT2}. The
situation is rather familiar: most of the difficulty lies in
guessing the loop group and a convenient representation of the dual
of its Lie algebra. In this case, the group is $LG^\Sigma_-$, the
set of smooth loops $g: {\Bbb S}^1 \to GL(N,\Bbb C)$ satisfying
the reality condition $\overline{g(z)}= g(\bar z)$ which admit
an analytic continuation to the exterior of the unit circle with
$g(\infty) = I,$ the identity matrix, and which are kept fixed by
the involution $\Sigma (g)(z) = (g(-z)^{T})^{-1}$. The nondegenerate
pairing on arbitrary loops
$$(X,Y)\, = \oint_{|z|=1} tr(X(z)Y(z))\,\frac{dz}{2\pi i}
                    = \sum_{j}\, tr (X_{j}Y_{-j-1})$$
identifies the Lie algebra dual $(L\frak g^{\sigma}_{-})^{*}$ with the
set of power series $ S_0  + N_1 z + \ldots$, where the
coefficients of the powers $z^k, k=0,1,\ldots$, are alternatively
symmetric and skew-symmetric real matrices.

The proof of Liouville integrability of the Bloch-Iserles equation
(*) takes two familiar forms. In one of them, embed $LG^\Sigma_-$
into the larger loop group $LG^\Sigma$ of all smooth invertible
loops and there search for a second bracket structure, induced by an
$R$-matrix, from which the existence of commuting integrals follows.
Otherwise, consider the Lie group $LG^\Sigma_{+} \times
LG^\Sigma_{-},$ corresponding to the Lie algebra anti-direct sum
$L\fg^{\sigma}_{+}\ominus L\fg^{\sigma}_{-}$  which as a manifold
can be identified with the subset of loops in $LG^\Sigma$ which
admit a Birkhoff factorization $g = g_{+}g_-^{-1}$ with trivial
diagonal part (\c{PS} is a good reference for the few facts about
Birkhoff factorization which will be used in this paper). The Lie
bracket of $L\fg^{\sigma}_{+}\ominus L\fg^{\sigma}_{-}$ is of course
identical to the bracket defined by the $R$-matrix above. Both
points of view will be described in an intertwined fashion in
Sections 2 and 3. Integrability is then proved in Section 4 by
showing that the integrals suggested by Bloch and Iserles are indeed
commuting, and generically span half of the dimension of the
coadjoint orbit of interest. A formula for the commuting flows in
terms of a Birkhoff factorization is the content of Section 5.

The loop group scenario embeds equation (*) for $n \times n$
matrices into a coadjoint orbit $\Cal{O}$ of generic dimension $2[
n^2/4]$ (here $[x]$ is the greatest integer less than or equal to
$x$). In Section 6, we consider a question raised by Bloch and
Iserles in their original paper: to find a finite dimensional matrix
group such that their equation is realized as a Hamiltonian system
on the dual of its Lie algebra. In \c{BI}, the authors presented an
interesting algorithm to study this problem, adapted from a
constructive proof of Ado's theorem, which states that any finite
dimensional Lie algebra admits a matricial faithful representation.
From the loop group setup, we obtain an explicit realization: there
is an $n^2$-dimensional nilpotent matrix group $G_f$, faithfully
realized on $3n \times 3n$ matrices, giving rise to a coadjoint
orbit isomorphic to $\Cal{O}$. The prescription yielding the finite
dimensional realization is simple and rather general. Since the
orbit $\Cal{O}$ is described in terms of (bi-)truncated Laurent
expansions of its elements, one should expect that the full group is
not needed to cover the orbit. Indeed, as we shall see, only the
first three terms of a loop in $LG^\Sigma_-$ are relevant to the
computation of the coadjoint action on the orbit of interest. This
suggests taking a quotient $G_f = LG^\Sigma_-/ LG^\Sigma_{-3}$ by
the normal subgroup $LG^\Sigma_{-3}$ of loops of the form $I +
O(z^{-3})$: it turns out that $G_f$ admits a simple Heisenberg-like
representation.

In the final section, we present an example of an infinite
dimensional orbit obtained by the same coadjoint action. Here, the
Hamiltonian which gives rise to equation (*) in the small orbit now
induces a system of partial differential equations containing some
simple integral terms.

Arieh Iserles informed us that a forthcoming paper, together with A.
Bloch, T. Ratiu and J. Marsden, will be dedicated to a different
proof of the integrability of the Bloch-Iserles equation, using very
different techniques.

The second author acknowledges support from CNPq and FAPERJ, Brazil.

\bigskip
\bigskip

\subhead
2. \ The Lie algebraic setup and the Poisson structure
\endsubhead
\bigskip

Let $G$ be $GL(N,\Bbb R)$ with Lie algebra $\fg$, and let $LG$ be the
group of loops $g:S^{1}\longrightarrow GL(N, \Bbb C)$ satisfying
the reality condition $\overline{g(z)} = g(\bar z)$.  We consider the
involutive automorphism
$\Sigma: LG\longrightarrow LG$,
given by $\Sigma(g)(z) = (g(-z)^{T})^{-1}$
and its fixed point set (the stable locus)
$$LG^{\Sigma} = \{\,g\in LG\mid g(z)(g(-z))^{T} = I\,\}.\eqno(2.1)$$
Clearly, $LG^{\Sigma}$ is a Lie subgroup of $LG.$  We shall denote by
$LG^{\Sigma}_{+}$ (resp.~ $LG^{\Sigma}_{-}$) the Lie subgroup of
$LG^{\Sigma}$ consisting of loops which extend analytically
to the interior (resp.~ to the exterior) of the unit circle,
with the additional requirement that loops in $LG^{\Sigma}_{-}$
take the value I, the identity matrix, at $\infty.$

Let $L\fg$ be the Lie algebra of $LG$, consisting of loops
$X(z) = \sum_{-\infty}^{\infty} X_{j} z^{j}$
with coefficients $X_{j}\in \fg.$
Then the Lie algebra $L\fg^{\sigma}$ of $LG^{\Sigma}$ is the stable locus
of the linearization $\sigma: L\fg \longrightarrow L\fg$ of
$\Sigma$, given by the formula
$$\sigma(X)(z) = -(X(-z))^{T} = \sum_{j} \theta(X_{j})(-z)^{j}\eqno(2.2)$$
where $\theta$ is the Cartan involution $\xi\mapsto -\xi^{T}.$
Therefore, if  $\fk$ (resp.~ $\frak p$) denote the $+1$ (resp.~
$-1$) eigenspace of $\theta$, consisting of skew-symmetric  (resp.~
symmetric) matrices, then explicitly,
$$L\fg^{\sigma} = \{\,X\in L\fg\mid X_{2j}\in \fk,\,X_{2j+1}\in \frak
p\,\,\,\, \hbox{for all}\,\, j\,\}.\eqno(2.3)$$

From the definition of $LG^{\Sigma}_{+}$ (resp.~ $LG^{\Sigma}_{-}$),
it is clear that its Lie algebra $L\fg^{\sigma}_{+}$
(resp.~$L\fg^{\sigma}_{-}$) consists of elements of the form
$\sum_{j\geq 0} X_{j}z^{j}$ (resp.~ $\sum_{j<0} X_{j}z^{j}$). Hence
the elements in $L\fg^{\sigma}_{-}$ are equal to $0$ at $\infty$ and
we have the splitting
$$L\fg^{\sigma} = L\fg^{\sigma}_{+} \oplus L\fg^{\sigma}_{-}\eqno(2.4)$$
with associated projection maps $\Pi_{+}$ and $\Pi_{-}$.

We now introduce
$$\TG = \lbrace g\in LG^{\Sigma}\mid g = g_{+} g^{-1}_{-},\,\hbox{where}\,\,
   g_{-}\in LG^{\Sigma}_{-},\, g_{+}\in LG^{\Sigma}_{+}
 \,\rbrace.\eqno(2.5)$$
Then an easy argument shows that  for $g\in \TG$, the factors
$g_{\pm}$ are unique (more details will be provided in Section 5).
Moreover, from the Birkhoff factorization theorem \c{PS}, $\TG$ is a
dense open subset of $LG^{\Sigma}$ in the natural topology.
Following the procedure in \c{DLT1},\c{DLT2}, we can endow $\TG$
with a Lie group structure by defining the multiplication
$$g\ast h \equiv g_{+}hg_{-}^{-1}.\eqno(2.6)$$
Clearly, the map $\Psi: \TG \longrightarrow LG^{\Sigma}_{+} \times
LG^{\Sigma}_{-}$ given by $\Psi(g) = \Psi(g_{+} g^{-1}_{-}) \mapsto
(g_{+},g_{-})$ is a Lie group isomorphism, when the image is
equipped with the product group structure. Consequently, the
pull-back of the standard Lie bracket on $L\fg^{\sigma}_{+} \oplus
L\fg^{\sigma}_{-}$ (Lie algebra direct sum) under
$\psi= T_{e}\Psi: X\mapsto(\Pi_{+}X, -\Pi_{-}X)$ yields the Lie bracket on
$\tg= Lie(\TG)$:
$$[X,Y]_{\tg} = [\,\Pi_{+}X, \Pi_{+}Y\,] - [\,\Pi_{-}X, \Pi_{-}Y\,]
\eqno(2.7)$$
for all $X,Y \in \tg.$
Similarly, one may use $\psi$ to obtain a formula for
the adjoint action of $\TG$ on $\tg$,
$$Ad_{\TG}(g) X= g_{-}(\Pi_{-}X)g^{-1}_{-}+ g_{+}(\Pi_{+}X)g^{-1}_{+} .\eqno(2.8)$$

\noindent{\bf Remark 2.1} By standard r-matrix theory \c{STS},
$$R = \Pi_{+} - \Pi_{-}\eqno(2.9)$$
is a solution of the modified Yang-Baxter equation, i.e.,
$$[RX, RY]-R([RX,Y] + [X,RY]) = -[X,Y]\eqno(2.10)$$
for all $X,Y\in L\fg^{\sigma}$.  By an easy computation, we
can show that the $R$-bracket on $ L\fg^{\sigma}$ given by
$$[X,Y]_{R} = {1\over 2}([RX,Y] + [X,RY])\eqno(2.11)$$
coincides with $[\cdot,\cdot]_{\tg}$.
Thus the vector space $L\fg^{\sigma}$ equipped
with the Lie bracket $[\cdot,\cdot]_{R}$ is identical to
the Lie algebra $\tg$.

\medskip

To identify the duals of the Lie algebras, we introduce
the following nondegenerate
invariant pairing on  $L\fg$:
$$\eqalign{(X,Y)\, &= \oint_{|z|=1} tr(X(z)Y(z))\,\frac{dz}{2\pi i}\cr
                    &= \sum_{j}\, tr (X_{j}Y_{-j-1}).\cr}\eqno(2.12)$$
As the reader will see, this choice of pairing is critical for what
we have in mind.  Using this pairing, we identify the
algebraic dual $(L\fg^{\sigma})^*$ of $L\fg^{\sigma}$  with the
stable locus of the Lie algebra anti-isomorphism $\sigma^*=-\sigma$:
$$L\fg^{\sigma^*} =\{\,X\in L\fg\mid X_{2j+1}\in \fk,\,X_{2j}\in \frak p
   \,\,\,\,\hbox{for all}\,\, j\,\}.\eqno(2.13)$$
Notice that the alternation of symmetric and skew-symmetric coefficients
in the Laurent expansion still holds, now with a parity opposite to the one
found in elements in the Lie algebra. By an easy computation making
use of (2.8) and (2.12), we find
$$Ad^{*}_{\TG}(g) A = \Pi_{-}(g^{-1}_{+}A g_{+}) +
\Pi_{+}(g^{-1}_{-}Ag_{-})\eqno(2.14)$$
where we have used the fact that $\Pi^{*}_{+} = \Pi_{-}$ and
$\Pi^{*}_{-} = \Pi_{+}.$
\smallskip
\noindent{\bf Remark 2.2} In the general context of an involutive
automorphism of a finite dimensional Lie algebra, the author in
\c{R} also considered the stable locus of the corresponding
involution on the loop algebra.  However, a different choice of
pairing was used in \c{R}, thus leading to a different identification
of the dual.

\medskip

Finally, the Lie-Poisson structure
for smooth functions on $\tilde{\fg}^*$ is given by the usual formula:
$$\{F_1,F_2\}(X) = (X, [dF_1(X), dF_2(X)]_{\tg}).\eqno(2.15)$$
\bigskip
\bigskip

\subhead 3. \ A Hamiltonian for the Bloch-Iserles equation
\endsubhead
\bigskip

The following result is standard.

\proclaim {Proposition 3.1} For $1\leq k\leq n$, $\ell\in \Bbb Z$,
define
$$H_{k\ell}(X) = {1\over (k+1)} \oint_{|z|=1} tr (X(z)^{k+1}) \ \frac{dz}{2\pi iz^{\ell+1}} ,
\eqno(3.1)$$ then the $H_{k\ell}$'s Poisson commute with respect to
$\{\cdot,\cdot\}.$ The Hamiltonian equation of motion generated by
$H_{k\ell}$ is given by
$$\dot X (z) = \left[\, \Pi_{+} ((X(z)^{k}z^{-(\ell+1)}),
X(z)\,\right].\eqno(3.2)$$
\endproclaim

\demo {Proof} Since $Ad^{*}_{g}(X)(z) = g(z)^{-1}X(z)g(z),$ the
Hamiltonians are invariant under the coadjoint action of
$LG^{\Sigma}$, i.e.
$$H_{k\ell}(Ad^{*}_{g}(X)) = H_{k\ell}(X),  \quad g\in LG^{\Sigma}. $$
By classical r-matrix theory, we then conclude that they Poisson commute.
The equation of motion for $H_{k\ell}$ is a straightforward
computation. \pf
\enddemo

Our next task is to search for interesting finite dimensional
coadjoint orbits in $\tg^*$. For any $m, n\in \Bbb Z_{+},$ define
$$\tg^{*}_{(m,n)} = \lbrace\, X\in L\fg^{\sigma^*}\mid X(z) = \sum_{j=-m}^{n}
    X_{j}z^{j} \,\rbrace.\eqno(3.3)$$
Clearly, we can extend this definition to the case where $m=0$
in which case we set $\tg^{*}_{n} = \tg^{*}_{(0,n)}$.

\proclaim
{Proposition 3.2} The sets $\tg^{*}_{(m,n)}$ are invariant
under $Ad^{*}_{\TG}(g)$ for any $g\in \TG.$
\endproclaim

\demo
{Proof} Take $X\in \tg^{*}_{(m,n)}$.  For $m =0$, the
analyticity property of $g_{\pm}$ implies that
$$\Pi_{-}(g^{-1}_{+}X g_{+}) = 0,$$
$$\Pi_{+}(g^{-1}_{-}X g_{-}) = X_{n}z^{n} + \sum_{j=0}^{n-1} B_{j}z^{j}$$
for some matrices $B_{j}$, and we are done (we have used $g_{-}(\infty) =I$
in deriving the second line above). For $m\in \Bbb Z_{+}$, note that
$$Ad^{*}_{\TG}(g)X = \Pi_{-}(g^{-1}_{+}(\Pi_{-}X) g_{+}) +
\Pi_{+}(g^{-1}_{-}(\Pi_{+}X) g_{-}).$$
Therefore, the same
calculation as before applied to $\Pi_{+}X$ shows that
$$\Pi_{+}(g^{-1}_{-}(\Pi_{+}X) g_{-}) = X_{n}z^{n} + \sum_{j=0}^{n-1}
B_{j}z^{j}$$
for some matrices $B_{j}$.  On the other hand, it is easy to
check that
$$ \Pi_{-}(g^{-1}_{+}(\Pi_{-}X) g_{+}) = g_{+}(0)^{T}X_{-m} g_{+}(0) +
\sum_{j= -m+1}^{-1}
    C_{j} z^{j}$$
for some matrices $C_{j}$.  Putting the above two expressions together,
the assertion follows.
\pf
\enddemo

As a consequence of this result, we obtain finite dimensional
coadjoint orbits through the elements in $\tg^{*}_{n}$ and
$\tg^{*}_{(m,n)}$ and the Hamiltonian equation in (3.2) restricts to
these orbits.  In the rest of the section, we shall focus on the
case $\tg^{*}_{1},$ in which the Bloch-Iserles equation lies.

\proclaim {Proposition 3.3} (a) Consider the loop $S_{0} + z
N_{0}\in \tg^{*}_{1}$. The $Ad^{*}_{\TG}$-orbit through $S_{0} +
zN_{0}$ is given by the affine linear space
$$
{\Cal O}_{S_{0} + zN_{0}} = \lbrace Ad^{*}_{\TG}(g) (S_{0}+
zN_{0})\mid g\in \TG\, \rbrace
   = \lbrace (S_{0} + [\,N_{0}, P]) + zN_{0}\mid P\in \frak
     p\,\rbrace.$$
     In particular, if $N_0$ has simple spectrum, the tangent space of
     ${\Cal O}_{S_{0} + zN_{0}}$
     is naturally identified with the vector space of real, symmetric
     matrices which are orthogonal to the matrix polynomials $p(N_{0}^2)$.
\newline
\noindent (b) The Hamiltonian equation of motion
$$(S + z N)^{.} = \left[\, \Pi_{+} ((S + zN)^{2} / z)), S + zN\,\right]
   \eqno(3.5)$$
generated by $H_{20}$ on the Poisson submanifold $\tg^{*}_{1}$ is
equivalent to the Bloch-Iserles equation
$$\eqalign{
           &\dot N = 0,\cr
           &\dot S = [\,NS + SN, S\,].\cr}\eqno(3.6)$$
\endproclaim

\demo {Proof}  We have
$$g_{-}(z) = I + g'_{-}(\infty) z^{-1} + O(z^{-2})$$
where $g'_{-}(\infty) = \frac {dg_{-}}{d z^{-1}}(z=\infty)$ and a
simple computation obtains
$${\Cal O}_{S_{0} + zN_{0}} = \lbrace S_{0} + [\,N_{0},
g'_{-}(\infty)] + zN_{0} \rbrace. $$ Since $g^{-1}_{-}(z)
=(g_{-}(-z))^{T},$ the matrix $g'_{-}(\infty)$ must be symmetric
and, conversely, any symmetric matrix $P\in \frak p$ is
$g'_{-}(\infty)$ for an appropriate loop $g_{-}$.  Now, endow $\frak
p$ with the usual inner product $<P,Q> = \ tr \ PQ$. Then the linear
map ${\Cal B}_{N_{0}}: \frak p \to \frak p$ taking $P$ to $[P,N_0]$
is skew symmetric and its range is orthogonal to its kernel, given
by the set of symmetric matrices commuting with $N_0$. If $N_0$ has
simple spectrum, the description of the tangent space to the orbit
then follows. The equation of motion for $H_{20}$ is easy to
compute, once one observes that only $S$ may vary in an orbit in
$\tg^{*}_{1}$. \pf
\enddemo
\smallskip
\noindent {\bf Remark 3.4} (a) From the proof of Proposition 3.2, it
should be clear that for $X\in\tg^{*}_{n}$, the equations in (3.2)
can be regarded as Hamiltonian systems on the coadjoint orbits of
the Lie group $LG^{\Sigma}_{-}$, when the dual
$(L\fg^{\sigma}_{-})^*$ of its Lie algebra is identified with the
collection of power series  $ S_0  + N_1 z + \ldots$, where the
coefficients of the powers $z^k, k=0,1,\ldots$ are alternatively
symmetric and skew-symmetric real matrices.  The advantage of using
the slightly more general formulation here is that it puts the
groups $LG^{\Sigma}_{+}$ and $LG^{\Sigma}_{-}$ on equal footing.
Otherwise, the relevance of $LG^{\Sigma}_{+}$ in the solution
formula in Proposition 5.1 may look somewhat mysterious.
\smallskip
\noindent (b) The Bloch-Iserles equation can be regarded as an isospectral
deformation of a general $n\times n$ matrix keeping the skew-symmetric part
fixed.  Indeed, if we set $z=1$ in
$(S + z N )^{.} = [\, NS + SN + z N^2 , S + z N\,]$
and let $M = S + N$, then a straightforward calculation shows
that
$$\dot M = {1\over 4} \left[\, (M^{T}M + MM^{T}) + (M^{T})^{2}, M\,\right].
\eqno(3.7)$$
Conversely, if $M$ satisfies the above equation, then
$$\dot M^{T} = {1\over 4} [\, (M^{T}M + MM^{T}) + M^{2}, M^{T}\,]\eqno(3.8)$$
and by a direct computation, we find $(M-M^{T})^{.} = 0$. Now, the
Bloch-Iserles equation preserves the spectral curve $\det (S + z N -
w) =0.$ (See Remark 4.3 below.)  On the other hand, the Toda flow in
\c{DLT1} also has a spectral curve, given by $\det (M + h(M^T -M)
-w) =0.$   If we write $M = S + N$, then by changing the variable
$h$ above in the obvious way, we can make the Toda curve coincide
with $\det (S + zN -w) =0$.  Percy Deift asked us if (3.7) might
arise from a combination of flows associated to the Toda curve. To
answer this question, note that any flow generated by a combination
of the coefficients from the Toda curve must be of the form $\dot M
= [\, \Pi_{\frak l}\, F(M, M^{T}), M \,].$ (See \c{DL}.) Here,
$\frak l$ is the Lie subalgebra of $\fg$ consisting of lower
triangular matrices and $\Pi_{\frak l}$ is the projection map to
$\frak l$ relative to the splitting $\fg = \frak k \oplus \frak l$,
and $F(M, M^T)$ is a polynomial in $M$ and $M^{T}$. If we denote the
right hand side of (3.7) by $X(M)$, then it is obvious that
$X(OMO^T) = OX(M)O^T$ for any orthogonal matrix $O.$ However, the
expression $[\, \Pi_{\frak l}\, F(M, M^{T}), M \,]$ clearly does not
satisfy this invariance property.  Hence the answer to the above
question is in the negative. From a different point of view, it is
also tempting to see if one can constrain the flows arising from the
Toda curve to the submanifold $Q = \{\,M \in \fg \mid M - M^T =N
\,\}$ where $N$ is a fixed skew-symmetric matrix. Unfortunately, the
submanifold $Q$ is not a cosymplectic submanifold in the sense of
Weinstein \c{W}. (The cosymplectic submanifolds are the
generalization in a Poisson context  where we can carry out Dirac's
idea of constraining a Hamiltonian system.) Thus this idea also
fails.

\bigskip
\bigskip

\subhead
4. \ Liouville integrability
\endsubhead
\bigskip

Let $sym_{ij}(A,B)$, the {\it $ij$-symmetrizer} of matrices $A$ and
$B$, denote the sum of all monomials in $A$ and $B$ of degree $i+j$
consisting of $i$ $A's$ and $j$ $B's$.

\proclaim {Lemma 4.1} (a) $[\ sym_{i,j+1}(A,B) \ , \ A \ ] + [\
sym_{i+1,j}(A,B)\ ,\ B \ ]=0$.
\newline
\noindent (b) Let $A=S$ be symmetric and $B=N$  skew-symmetric. Then
$sym_{ij}(S,N)$ is symmetric (resp. skew-symmetric) if $j$ is even
(resp. odd).
\newline
\noindent (c) If $A$ and $B$ are $n \times n$ matrices, than each
symmetrizer $sym_{n - \ell,\ell}(A,B),\ell=0,\ldots,n,$ is a linear
combination of symmetrizers $sym_{r,s}(A,B)$, of smaller degree $r+s
< n$.
\newline
\noindent (d) For generic choices of $A$ and $B$ (i.e., for an open,
dense set of pairs $(A,B)$), $sym_{k-\ell,\ell}(A,B),
k,\ell=0,\ldots,n-1,$ are linearly independent. The same is true for
generic choices of $A=S$ symmetric and $B=N$ skew-symmetric.
\endproclaim

\demo {Proof} To simplify notation, we omit the obvious matrix
dependence of $sym_{ij}$. Statement (a) follows from
$$ sym_{i+1,j+1} = A  \ sym_{i,j+1} + B \ sym_{i+1,j} = sym_{i,j+1} \  A+
sym_{i+1,j} \  B.$$ The first equality states that any monomial in
$sym_{i+1,j+1}$ either starts with $A$ and is completed with $i \
A's$ and $j+1 \ B's$ or starts with $B$ and is completed with $i+1 \
A's$ and $j \ B's$. The second equality is obtained by taking into
account the last matrix of each monomial. The proof of (b) is
obvious. To see (c), notice that the claim for $i=0$ (resp. $i=n$)
is just the Cayley-Hamilton theorem for $B$ (resp. $A$). More
generally, write the Cayley-Hamilton theorem for the matrix $M = A +
z B$,
$$ - M^n = \sigma_1 M^{n-1}+ \sigma_2 M^{n-2}+ \ldots + \sigma_n M^{0},$$
where  $\sigma_j = \sigma_j(z)$, a standard symmetric function of
the eigenvalues of $M$, is a polynomial of degree $j$ in $z$. Now
collect terms in $z$: on the left hand side, the coefficient of
$z^j$ is of the form $sym_{n-\ell,\ell}$; on the right hand side, it
is a linear combination of symmetrizers of degree smaller than $n$.
To prove (d), let $B_*$ be the matrix whose only nonzero entries,
equal to one, lie along the  subdiagonal of entries with indices
$(r,s)$ satisfying $r-s=1$ and let $A_*$ be a diagonal matrix with
entries along the  diagonal in geometric progression
$c^1,c^2,\ldots, c^n$. Clearly, $sym_{k-\ell, \ell}$ only has
nonzero entries along the subdiagonal $r-s = \ell$, forming a
geometric progression of order $c^{k-\ell}$. Thus, to check if a
linear combination of symmetrizers is independent, it suffices to
consider the independence of subsets of symmetrizers consisting of a
fixed index $\ell$ (keep in mind that, by hypothesis, $\ell < n $),
for which $k-\ell$ runs from 0 to $n-\ell-1$. We thus have to prove
the linear independence of $n-\ell$ vectors in ${\Bbb R}^{n-\ell}$,
with coordinates in geometric progressions of ratio $c^1,\ldots,
c^{n-\ell-1}$: this is clearly true for a generic choice of $c$ (we
thank Nicolau Saldanha for suggesting matrices $A_*$ and $B_*$).
From analyticity, independence holds for generic choices of $A$ and
$B$, proving the first part of (d). We now show that generic
independence still holds for symmetric (resp. skew-symmetric)
choices of $A$ (resp. $B$). By continuity, independence still holds
for matrices sufficiently close to $A_*$ and $B_*$. In particular,
this is still true if we keep $S_* = A_*$ and change $B_*$ by adding
a small negative number on its subdiagonal of entries $(r,s)$ for
which $r-s = -1$, giving rise to a matrix $\tilde {B_*}$. Now,
consider a diagonal matrix $D$ so that $N_* = D {\tilde {B_*}}
D^{-1}$ is skew-symmetric. The pair $S_* = D A_* D^{-1} $ and $N_*$
also has independent symmetrizers of degree less than $n$ and
genericity for pairs $(S,N)$ follows again by analiticity.\pf
\enddemo

As usual, let $S$ and $N$ be symmetric and skew-symmetric $n \times
n$ matrices. The set of symmetrizers $sym_{k-\ell,\ell}(S,N)$ having
degree $k < n$ which are symmetric matrices (i.e., those for which
$j$ is even, by the previous lemma) splits in two sets. The first
will give rise to the family of commuting integrals $H_{k\ell}$, the
second relates to the Casimirs of generic coadjoint orbits within
$\tg^{*}_{1}$:
$$\eqalign{
{\Cal I} &= \lbrace sym_{k-\ell,\ell}, \   1 \leq k \leq n-1, \ 0
\leq \ell \leq n-2,\  \ell \hbox{ even}\rbrace, \cr {\Cal S}_{N}
&=\lbrace sym_{\ell\ell},\ 0 \leq \ell \leq n-1, \ \ell \hbox{ even}
\rbrace. \cr}$$

From Lemma 4.1, for a generic choice of $S$ and $N$, the map ${\Cal
B}_N$ is injective when restricted to the vector space spanned by
${\Cal I}$: in particular, the matrices $[ \ sym_{k-\ell,\ell}(S,N)
\ , \ N \ ]$, for $\ sym_{k-\ell,\ell} \in {\Cal I},$ are linearly
independent.

\proclaim
{Theorem 4.2} (a) Let $N_0$ have simple spectrum and consider
an element $S_{0} + z N_{0}\in \tg^{*}_{1}$. Then ${\Cal O}_{S_{0} +
z N_{0}}$ is a $2[n^{2}/4]$ dimensional affine linear space given by
$$\left\{ (S + z N)\in \tg^{*}_{1}\mid N = N_{0},\,tr (SN^{\ell})=
tr (S_{0}N_{0}^{\ell}), \ell \hbox{ even },  0\leq\ell \leq
n-1\,\right\}. \eqno(4.2)$$
\newline
\noindent (b) The equation of motion of the Hamiltonian $H_{k\ell}$
is given by
$$ \dot{S} = [ \ sym_{k - (\ell+1),\ell+1}(S,N) \ , \  S  \ ]
= - [ \  sym_{k - \ell,\ell}(S,N) \ , \ N \  ], \quad \dot{N}=0.$$
The Hamiltonians in the set
$$\lbrace H_{k\ell},\  1 \leq k \leq n-1, \
0 \leq \ell \leq n-2,\  \ell \hbox{ even}\, \rbrace$$ are
generically independent on ${\Cal O}_{S_{0} + z N_{0}}$ and provide
$[n^2/4]$ commuting integrals.
\endproclaim

\demo {Proof} Statement (a) follows from the characterization of the
tangent space to the orbit given in Proposition 3.3, combined with
the independence of the elements of ${\Cal C}$. The first formula
for the equation of motion for $H_{k\ell}$ follows by expanding
$X(z)^{k}$ in the equation of motion as given in Proposition 3.1.
The second formula is a consequence of Lemma 4.1(a). Finally, by
considering the second description, generic independence of the
vector fields is a consequence of the remarks just above the
statement of the theorem. \pf
\enddemo

\noindent {\bf Remark 4.3} Equivalently, we can consider the
conserved quantities given by
the coefficients of the characteristic polynomial $p(z,w) = det(S +
zN - wI)$. Notice that  $(S + zN)^{T} = S -z N$, so $p(z,w)$ is an
even function of $z:$
$$p(z,w) = det(S + zN -wI) = \sum_{r=0}^{N}\sum_{k=0}^{\left[r\over 2\right]}
I_{rk}(S,N)z^{2k} w^{N-r}.$$ The fact that the functions $I_{rk}$
Poisson commute also follows by coadjoint invariance. Generic
independence follows from the generic independence of the
$H_{k\ell}$.  Note that the functions $I_{2k,k}$, $k= 1,\ldots,
[n/2]$ clearly depends only on $N$ and therefore are trivial
integrals.  On the other hand, it is easy to see that $\{I_{2k+1,
k}\}_{0\leq k\leq \left[{\frac{n+1}{2}}\right]-1}$ is equivalent to
the coadjoint orbit invariants $\{tr(SN^{2k})\}_{0\leq k \leq
[\frac{n-1}{2}]}$.  Hence the nontrivial integrals are given by
$\{I_{rk}\}_{0\leq k\leq [r/2]-1, 1\leq r\leq n}$ and this also
gives a total of
$$\sum_{r=1}^{n} \left[\frac{r}{2}\right] = \left[\frac{n^{2}}{4}\right]$$
conserved quantities, as required.
The fact that the {\it spectral curve} $p(z,w)=0$ is preserved
by the Bloch-Iserles equation means that the corresponding
flow is linearized on the corresponding Jacobian variety
and that the solution can be explicitly written down in
terms of Riemann theta functions.  We shall leave the details
to the interested reader. (See, however, Proposition 5.2 and Remark 5.3
in this connection.)

\bigskip
\bigskip

\subhead
5. \ Solution by factorization
\endsubhead
\bigskip

The flows associated to the integrals $H_{k\ell}$ may be described
by an explicit formula. Let ${\Bbb D}_+, {\Bbb S}^1$ and ${\Bbb
D}_-$ be respectively, the sets of complex numbers $z$ together with
$z = \infty$ for which $|z| \le 1$, $|z| = 1$ and $|z| \ge 1$.
Recall that (\c{PS}) a loop $\gamma \in LG$ admits a {\it Birkhoff
factorization} $\gamma(z) = \gamma_+(z) \ d(z) \ \gamma_-^{-1}(z).$
Here, $\gamma_+$ and $\gamma_-$ are restrictions to ${\Bbb S}^1$ of
analytic functions extending to the boundary of ${\Bbb D}_+$ and
${\Bbb D}_-$ taking their values on $GL(N,\Bbb C)$, the matrix
$d(z)$ is diagonal with diagonal entries of the form $z^{a_i}, a_i
\in {\Bbb Z}$, and $\gamma_-(\infty) = I$, the identity matrix. The
only ingredient which is not automatic in the derivation of the
formula for the flows is the proof of the triviality of the diagonal
factor in the factorization of the loops associated to the initial
conditions.

As stated, the Birkhoff factorization of an invertible loop is not
necessarily unique. What is true, and follows from an argument
similar to the given in the proof of the proposition below, is that
different factorizations have the same diagonal factor, up to
permutation of its diagonal entries. The loops we have to consider,
however, admit a special symmetry. We follow in spirit the arguments
in \c{GK}.

\proclaim {Proposition 5.1} Let $\gamma$ be a loop satisfying the
symmetry $\gamma(z)\gamma^T(-z)=I$.  Then its Birkhoff factorization has
trivial diagonal factor, $d(z) = I$, and is unique. Moreover, the nontrivial
factors also satisfy the symmetry,
$\gamma_{\pm}(z)\gamma_{\pm}^T(-z)=I$.
\endproclaim

\demo {Proof} Suppose $\gamma(z) = \gamma_+(z) d(z)
\gamma_-^{-1}(z)$, with the notation above: we have to show that all
the integers $a_i$ are equal to 0. For a continuous function $f:
{\Bbb S}^1 \to \ {\Bbb C}^*$, define its winding number $w(f)$ to be
the (signed) number of turns of its image around the origin. The
winding number is invariant under continuous deformations $f_t$
through functions which avoid the origin and is additive with
respect to products, $w(fg) = w(f) + w(g)$. Since $\gamma_+$ (resp.
$\gamma_-$) extends to ${\Bbb D}_+$ (resp. ${\Bbb D}_-$), taking
values at invertible matrices, one may deform the loop $\gamma$ to a
constant loop through loops taking values on invertible matrices, so
that $w(\det \gamma_+) = w(\det \gamma_-) = 0$. By additivity,
$w(\det \gamma) = w(\det\, d) = \sum_i a_i.$ On the other hand,
since $\gamma(z)\gamma^T(-z)=I$ and $w(\det \gamma(z))= w(\det
\gamma^T(-z)),$ we must have $w(\det \gamma) = 0$. Thus, $\sum_i
a_i= 0.$ Also from $\gamma(z)\gamma^T(-z)=I$, we have
$$ d(z) \gamma_-^{-1}(z)
\gamma_-^{-T}(-z) d(-z) = \gamma_+^{-1}(z) \gamma_+^{-T}(-z).$$
Equating diagonal entries of both sides, we obtain, for
$i=1,\ldots,n,$
$$   (\gamma_-^{-1}(z) \gamma_-^{-T}(-z))_{ii}
= (-1)^{a_i}z^{-2{a_i}}(\gamma_+^{-1}(z) \gamma_+^{-T}(-z))_{i,i},$$
which we denote, with the obvious attributions, by $ f_-(z)  =
z^{-2{a_i}}g_+(z)$. Clearly, the hypothesis $\gamma_-({\infty})=I$
implies $f_-(\infty)= 1$. Suppose $a_i < 0$ for some $i$. Then the
function $h(z)$, which agrees with $f_-$ on ${\Bbb D}_-$ and with
$z^{-2{a_i}} f_+(z) $ on ${\Bbb D}_+$ is an entire, bounded function
--- a constant --- satisfying $h(0)=0$ and $h(\infty)= 1$.
The upshot is that $a_i \ge 0$, and since $\sum_i a_i = 0,$ we must
have $a_i = 0$ for all $i$. A similar argument involving Liouville's
theorem obtains uniqueness of factorization. Symmetry of
$\gamma_{\pm}$ in turn follows from unique factorization applied to
the equation $\gamma(z)= \gamma_+(z) \gamma_-^{-1}(z) =
\gamma(-z)^{-T} = \gamma_+(-z)^{-T} \gamma_-(-z)^{T}.$ \pf
\enddemo

As we noted earlier, the vector space $L\fg^{\sigma}$ equipped with
the Lie bracket in (2.11) coincides with $\widetilde \fg$. Since the
modified Yang-Baxter equation (eqn. (2.10)) is a factorization
condition, our next result is standard (\c{STS}). We sketch the
argument in order to check that the symmetries required in the
previous proposition to obtain the Birkhoff factorization with
trivial diagonal factor indeed hold.

\proclaim {Proposition 5.2} Let  $f_{k\ell}(x,z)=
x^{k}{z^{-(\ell+1)}}$, for $\ell$ even. Then the solution of the
equation of motion
$$\dot X (z) = \left[\, \Pi_{+} ((X(z)^{k}z^{-(\ell+1)}),
X(z)\,\right]= \left[\, \Pi_{+} f_{k\ell}(X(z),z),X(z)\,\right],
\quad X(0,z) = X_0(z),$$ generated by the Hamiltonian $H_{k\ell}$ is
given by
$$ X(t,z)
= g_-^{-1}(t,z) X_0(z) g_-(t,z)= g_+^{-1}(t,z) X_0(z) g_+(t,z),$$
where $g_+$ and $g_-$ are obtained from the Birkhoff factorization
$$exp(-t f_{kl}(X_0(z),z)) = g_+(t,z) g_-^{-1}(t,z).$$
\endproclaim

\demo {Proof} First, notice that, for any real $t$, the loop
$\gamma(t,z) = exp(-t f_{k\ell}(X_{0}(z),z))$ satisfies the
hypothesis of the previous proposition. Indeed, it suffices to check
that, for even $\ell$, $\delta(z) = -f_{k\ell}(X_{0}(z),z)$
satisfies the linearization $\delta(z) + \delta^T(-z) = 0$, which is
obvious. Thus $g_-$ and $g_+$ are uniquely determined from
$exp(-tf_{k\ell}(X_0(z),z)) = g_+(t,z) g_-^{-1}(t,z)$.
Differentiating the above expression, we obtain
$$\aligned
&exp(-t f_{k\ell}(X_0(z),z)) \
f_{k\ell}(X_0(z),z) \\
= &\ \dot g_+(t,z) g_-^{-1}(t,z) - g_+(t,z)g_-^{-1}(t,z) \dot
g_-(t,z)g_-^{-1}(t,z) \endaligned
$$
so that
$$g_+^{-1}(t,z)\dot g_{+}(t,z) = -\Pi_+ (g_+^{-1} (t,z)f_{k\ell}(X_0(z),z)
g_+(t,z)) = -\Pi_+ f_{k\ell}(X(t,z),z)$$ for  $ X(t,z) =
g_+^{-1}(t,z)X_0(z) g_+(t,z)= g_-^{-1}(t,z)X_0(z) g_-(t,z)$. Taking
the derivative of the first expression for $X(t,z)$, we obtain
$$ \dot X(t,z) = \left[ \,  -g_+^{-1}(t,z)\dot g_+(t,z), X(t,z) \,
\right] = \left[\, \Pi_{+} f_{k\ell}(X(t,z),z)), X(t,z)\,\right].$$ \pf
\enddemo
\smallskip
\noindent{\bf Remark 5.3} The Birkhoff factorization for
$g_{\pm}(t,z)$ above can be solved explicitly in terms of Riemann
theta functions.  See, for example, \c{RS} and \c{DL} for details.
\smallskip
\noindent{\bf Remark 5.4} In the nonperiodic Toda hierarchy, there
are vector fields whose solution for a given initial condition $S_0$
are especially simple at integer times. One example is associated to
the $QR$ factorization of $exp(t f(S_0))$ for $f(x) = \ln (x)$. This
evolution is related to the $QR$ iteration, a starting point of many
algorithms to compute eigenvalues (\c{S},\c{DNT}). In a similar
fashion, there are evolutions related to some of the Hamiltonians
$H_{k\ell}$ which are algebraically solvable at integer times, since
the Birkhoff factorization  may be performed explicitly on loops of
matrices with polynomial entries.

\bigskip
\bigskip

\subhead 6. \ A finite dimensional group for the Bloch-Iserles
equation
\endsubhead
\bigskip

We now indicate how to obtain a finite dimensional group which
induces a coadjoint orbit diffeomorphic to ${\Cal O}_{S_{0} +
zN_{0}}$, on which the Bloch-Iserles equation arises as a
Hamiltonian system. This construction answers a question addressed
in \c{BI}.

Represent a loop $g(z)=(I + \frac{g_{-1}}{z} + \frac{g_{-2}}{z^2} +
\ldots) \in LG^\Sigma_-$ (set $g_0 = I$) as the bi-infinite
(convolution) matrix $M_g$ which, in Fourier variables, corresponds
to multiplying a (matrix) function defined on the circle by $g(z)$.
The matrix $M_g$ is obtained by prescribing the same basis $\lbrace
\ldots, z^{-2},z^{-1},z^0,z,z^2,\ldots \rbrace$ on domain and range.
It is block upper triangular, with $n \times n$ blocks indexed by
$(M_g)_{ij}$, equal to $g_{-k}$ for $i-j=-k$. Upper unipotent
convolution matrices form a group $\Cal M$, and matrices associated
to loops admitting the symmetry $g(z) g^T(-z)=I$ form a subgroup
${\Cal M}^\Sigma$. Consider the normal subgroup ${\Cal M}_3^\Sigma$
of ${\Cal M}^\Sigma$ on which the diagonals associated to $k=1$ and
$2$ have zero entries: such matrices correspond to loops of the form
$I + O(z^{-3})$. The quotient ${\Cal M}^\Sigma/{\Cal M}_3^\Sigma$
corresponds to 'forgetting' diagonals associated to $k \ge 3$ and is
clearly isomorphic to the group $G_f = LG^\Sigma_-/ LG^\Sigma_{-3}$
defined in the Introduction. Elements in the quotient may be
represented by matrices $M_g^0$, on which diagonals for which $k \ge
3$ are set equal to zero. Now, set $\Pi_3$ denote the orthogonal
projection on the three basis elements $1,z,z^2$. A simple
computation shows that ${\Cal M}^\Sigma/{\Cal M}_3^\Sigma$ is
(group) isomorphic to the $3n \times 3n $ matrices of the form
$\Pi_3 M_g^0 \Pi_3^*$ --- this is the finite group which will induce
the Bloch-Iserles equation, as we shall see.

We now provide the concrete description of ${\Cal M}^\Sigma/{\Cal
M}_3^\Sigma$. Consider the group of real $3n \times 3n$ matrices of
the form
$$g= g(S,N) = \pmatrix
I & S & \frac{S^2}{2} + N \cr 0 & I & S \cr 0 & 0 & I \cr
\endpmatrix ,$$
where the entries are $n \times n $ matrices. Here, matrices
$S,T,\tilde S, \tilde T$ are symmetric, $N,M,\tilde N, \tilde M$ are
skew symmetric. A simple computation gives $g(S,N) g(T,M) = g(S+T,
N+M + \frac{[S,T]}{2}),$ in agreement with the loop product
$$\eqalign{(I + &\frac{S}{z} + \frac{{\frac{S^2}{2}}+N}{z^2} + \ldots)
(I + \frac{T}{z} + \frac{{\frac{T^2}{2}}+M}{z^2} + \ldots)= \cr &I +
\frac{S+T}{z} + \frac{{\frac{(S+T)^2}{2}}+N+M +\frac{[S,T]}{2}}{z^2}
+ \ldots. \cr}$$ The Lie algebra $\fg_f$ and its dual $\fg_f^*$
consist of elements of the form
$$X(S,N) = \pmatrix
0 & S &  N \cr 0 & 0 & S \cr 0 & 0 & 0 \cr
\endpmatrix  \hbox{ and }
A(S,N) = \pmatrix 0 & 0 &  0 \cr \tilde S & 0 & 0 \cr 2 \tilde N &
\tilde S & 0 \cr
\endpmatrix.$$
Here the nondegenerate pairing is $(X, A) = tr X A$ and the usual
computations yield $Ad^*_{g(T,M)}(A(S,N)) = A(S + [N, T], N).$ Up to
trivial identifications, this is exactly the formula for the
coadjoint action of  $LG^{\Sigma}_{-}$ through the loop $S + zN$ described
in Section 3.
Clearly, the same Hamiltonians defined in Section 3 induce the same
flows in this setup.

\bigskip
\bigskip

\subhead 7. \ Another coadjoint orbit
\endsubhead
\bigskip

From formula (2.18), it is clear that an $n \times n$ matrix $N_0$
with a large kernel gives rise to especially small coadjoint orbits.
or a concrete example, split an arbitrary symmetric matrix $S$ into
blocks,
$$S= \pmatrix
a & b & u^{T} \cr
b & c & v^{T} \cr
u & v & B \cr
\endpmatrix \eqno(6.1)$$
where $a$, $b$ and $c$ are real numbers, $u$ and $v$ are vectors of
dimension $n-2$ and $B$ is a real, symmetric matrix of dimension
$n-2$. A direct computation shows that the flow $\dot S= [N,S^2]$ is
equivalent to the equations
$$\eqalign{\dot a \, &= 2 \langle u,  v
\rangle, \cr
          \dot  b \, &=  \langle v,  v \rangle - \langle u,  u \rangle, \cr
           \dot c \, &= - 2 \langle u,  v \rangle = - \dot a, \cr
            \dot B\, &= 0.\cr}\eqno(6.2)$$
With the same notation as above, the equation $\dot S= [N, S^{3}]$ reads
$$\eqalign{\dot a \, &= 2 b^3 + 2 b a^2
+2 b (\langle v ,v \rangle + \langle u ,u \rangle) + 2\langle u,Bv \rangle, \cr
\dot b\, &= -2a^3-2ab^2-
-2 a (\langle u ,u \rangle + \langle v ,v \rangle) + \langle v,Bv \rangle -
\langle u,Bu \rangle, \cr
\dot c\, &= -\dot a, \cr
\dot u\, &= (b^2 + a^2) v + b Bu -a Bv
+ \langle u,v \rangle u + \langle v,v \rangle v + B^2 v, \cr
\dot v\, &= -(b^2 + a^2)u - a Bu -b Bv
- \langle u,u \rangle u - \langle u,v \rangle v - B^2 u.
\cr}\eqno(6.3)$$

We select a simple special case. Suppose $u$ and $v$ are points in a
possibly complex vector space $V$. Say $B$ is a symmetric matrix
(i.e., $B^{T} = B$) acting on $V$, so that its entries may be
complex, and suppose that $V$ splits in two subspaces $V_e$ and
$V_o$ (we will denote them by even and odd vectors), which are real
orthogonal to each other and are interchanged by $B$. In particular,
we have $\langle u, B^{2k+1} v \rangle = 0$, for all natural $k$,
for $u$ and $v$ with the same parity. Finally, suppose that, at time
0, the entries $a, b$ and $c$ are 0, and the vectors $u$ and $v$
have the same parity. Then we must have that $a, b$ and $c$ are kept
constant equal to $0$ and $u$ and $v$ preserve their parity. Indeed,
the requirement that $a, b$ and $c = -a$ are constant equal to zero
is compatible with their evolutions, and the remaining equations,
keeping the notation above to denote the real inner product on $V$,
become
$$\eqalign{
\dot u\, &= \langle u,v \rangle u + \langle v,v \rangle v + B^2 v, \cr
\dot v\, &= \langle u,u \rangle u - \langle u,v \rangle v - B^2 u, \cr}
\eqno(6.4)$$
which respect parity. From uniqueness of solutions for a differential equation,
whatever solves this smaller system actually is the unique solution of the
original system provided $a, b$ and $c$ stay put.

We consider an infinite dimensional version of the above system.
Take $B = i D_x$, so that $B^2 = - D_{xx}$, and take $u(x)$ and
$v(x)$ to be simultaneously either even or odd functions in the
variable $x$. Let the inner product denote the usual real product in
$L^2(dx)$. The evolution equations for the functions $u(t,x)$ and
$v(t,x)$ read
$$\eqalign{
u_t \, &= \langle u,v \rangle u + \langle v,v \rangle v - v_{xx}, \cr
v_t \, &= \langle u,u \rangle u - \langle u,v \rangle v + u_{xx}. \cr}
\eqno(6.5)$$
In particular, this integro-differential system preserves the original parity
and the reality of the initial conditions.

Other partial differential equations which fit the format $\dot
S(t)=[N, S^k]$ may be obtained as follows. Set $S = M_q$ the
operator given by multiplication by the function $q(x)$. Choose now
$N = D M_\alpha + M_\alpha D$, where $D$ is the partial derivative
with respect to $x$. Clearly, $N$ is skew-symmetric and the
evolution equation for the operators becomes $q_t = 2 \alpha
(q^k)_x.$ These differential equations admit hierarchies of
infinitely many conserved quantities, but some of these conserved
quantities generate vector fields which are not differential
equations: in a sense, the coadjoint orbit through the initial
condition is too large.

\enddocument